\documentclass[10pt]{article}   	% use "amsart" instead of "article" for AMSLaTeX format
\usepackage{geometry}                		% See geometry.pdf to learn the layout options. There are lots.
\usepackage{graphicx}				% Use pdf, png, jpg, or eps§ with pdflatex; use eps in DVI mode
	\usepackage{mathrsfs}							% TeX will automatically convert eps --> pdf in pdflatex		
\usepackage{amssymb}
\usepackage{hyperref}
\usepackage{amsmath}
\usepackage[all,cmtip]{xy}
\usepackage{amscd,amssymb,amsfonts,amsmath,latexsym,amsthm}
\usepackage[nameinlink, capitalise]{cleveref}
\usepackage{etoolbox}
\title{Coarse geometry}
\author{Ulrich Bunke\thanks{Fakult{\"a}t f{\"u}r Mathematik,
Universit{\"a}t Regensburg,
93040 Regensburg,
ulrich.bunke@mathematik.uni-regensburg.de}
}
\newcommand{\Fun}{\mathbf{Fun}}

\newcommand{\Top}{\mathbf{Top}}
\newcommand{\Set}{\mathbf{Set}}
\newcommand{\Unif}{\mathbf{Unif}}
\newcommand{\Coarse}{\mathbf{Coarse}}
\newcommand{\BC}{\mathbf{BornCoarse}}
\newcommand{\cC}{\mathcal{C}}
\newcommand{\cB}{\mathcal{B}}
\newcommand{\cY}{\mathcal{Y}}
\newcommand{\diag}{\mathrm{diag}}
\newcommand{\R}{\mathbb{R}}
\newcommand{\Z}{\mathbb{Z}}
\newcommand{\C}{\mathbb{C}}
\newcommand{\nat}{\mathbb{N}}
\newcommand{\Orb}{\mathrm{Orb}}
\newcommand{\bM}{\mathbf{M}}
\newcommand{\pr}{\mathrm{pr}}
\newcommand{\colim}{\mathrm{colim}}
\newcommand{\map}{\mathrm{map}}
\newcommand{\cX}{\mathcal{X}}
\newcommand{\Mod}{\mathbf{Mod}}
\newcommand{\bV}{\mathbf{V}}
\newcommand{\Ccat}{C^{*}\mathbf{Cat} }
\newcommand{\Sp}{\mathbf{Sp}}
\newcommand{\End}{\mathrm{End}}
\newcommand{\Fin}{{\mathcal{F}in}}
\newcommand{\Dirac}{ D\hspace{-0.24cm}/ \hspace{0.08cm}}
\newcommand{\ind}{\mathrm{index}}
\newcommand{\free}{\mathrm{free}}
\newcommand{\Hilb}{\mathbf{Hilb}}
\newcommand{\hB}{\hfill $\Box$}

\newcommand{\id}{\mathrm{id}}
\newcommand{\Var}{\mathrm{Var}}

\newcommand{\coarse}{{\tiny c}}

\newtheorem{theorem}{Theorem}[section] 
\newtheorem{prop}[theorem]{Proposition}

\newtheorem{ddd}[theorem]{Definition}

\theoremstyle{remark}
\theoremstyle{definition}

\newtheorem{ex}[theorem]{Example}
\newtheorem{rem}[theorem]{Remark}

\begin{document}
\maketitle
\begin{abstract}
In this chapter we introduce coarse and bornological coarse spaces. We explain the concept of a coarse homology theory and discuss the examples of coarse ordinary homology and coarse $K$-homology in some detail.
\end{abstract}

\maketitle\tableofcontents
\setcounter{tocdepth}{5}

%\section{}\cite{higson_pedersen_roe}
%\subsection{}

\section{Introduction}
Coarse geometry was invented by J. Roe \cite{MR1147350}, \cite{roe_lectures_coarse_geometry}. The original motivation came from  index theory of Dirac type operators on complete Riemannian manifolds \cite{roe_index_coarse}. But coarse geometry  is also a framework to study geometric properties of groups and assembly maps for algebraic $K$-theory  \cite{Bartels_2004},   \cite{blr}. 
Via the cone  construction  \cite{higson_pedersen_roe}, \cite{MR1834777} it subsumes the controlled topology approach   
 \cite{Carlsson_1995}, \cite{Weiss_2002} to algebraic $K$-theory. Coarse geometry also has been used to study topological insulators in mathematical physics \cite{Ewert_2019}, \cite{Ludewig_2021}.   
%\cite{MR1834777}, 
%\cite{Carlsson_1995}, \cite{Weiss_2002}

\begin{itemize}
\item In this survey we explain the category of $G$-bornological coarse spaces \cite{equicoarse}
 as a basic framework  to study large scale invariants of metric spaces with $G$-action. 
 \item 
 We introduce the concept of a coarse equivalence 
 and  the idea of a coarse invariant.  As examples we discuss the sets of coarse components and the Higson corona.
 \item We introduce the concept of an equivariant  coarse homology theory.
 
\item We  provide a complete description  of  equivariant coarse ordinary homology and equivariant coarse $K$-theory.
   \end{itemize}
  
  Let us point out that
though coarse homology theories constitute a central aspect of coarse geometry  this survey only touches a very small portion of the field of coarse geometry. Furthermore, the list of given references 
is by far not complete and just provides entry points for further reading.

\section{Coarse spaces}						% Activate to display a given date or no date

In metric geometry, the mutual relation between points of a set is encoded by a distance function $d:X\times X\to [0,\infty]$, where the generalised real  number $d(x,y)$ is interpreted as the distance from $x$ to $y$. %Thereby  one requires that: \begin{enumerate}  \item $d(x,y)=0$ if and only if $x=y$, \item  the symmetry $d(x,y)=d(y,x)$, \item and the triangle inequality 
%$d(x,z)\le d(x,y)+d(y,z)$. \end{enumerate}
The pair $(X,d)$ is called a metric space. Note that in contrast to the usual conventions we allow infinite distances in order to model spaces with more than one coarse component later on.

The space $\R^{n}$ with the euclidean distance $d_{eu}$ is an example of a metric space. Every subset of a metric space becomes a metric space with the restricted distance function.

A distance function  often encodes much more   structure about the geometry of $X$ as one is interested in or can determine in practise. 

\begin{ex}Assume that $G$ is a group. If we choose a subset $S$ of $G$, then we can define a distance $d_{S}(g,h)$ as the minimal number of elements of $S\cup S^{-1}$ needed to express $gh^{-1}$. This could be infinite if $S$ does not generate $G$. But even if $S$ and $T$ are finite generating sets of $G$, then the distance functions $d_{T}$ and $d_{S}$ are different in general. But one is  only interested  in geometric properties of $G$ that do 
  not depend on that choice. \hB \end{ex}

Let $(X,d)$ be a metric space.
If one wants to concentrate  on the small scale structure of $X$ locally, then one usually only considers the topology determined by the distance function which is generated by the open balls
\begin{equation}\label{vwewervijovwvwervwerv}B(x,r):=\{y\in X\mid d(x,y)<r\}
\end{equation}  for all $r$ in $(0,\infty)$ and $x$ in $X$. If one wants to be able to compare the local scales at different points of $X$, then it is natural to work with the uniform structure on $X$ generated by the entourages
 \begin{equation}\label{sbgbsgbvfsfvfvsfdv}U_{r}:=\{(x,y)\in X\times X\mid d(x,y)<r\}
\end{equation}  for  all $r$ in $(0,\infty)$.   Similarly, 
in coarse geometry one is  interested in the large scale structure of the metric space only.
 
A mathematical way to introduce structures is to describe a category whose objects represent the structures of interest, and whose morphisms are structure preserving maps. In the cases discussed above we arrive at the categories $\Top$ of topological spaces and continuous maps and $\Unif$ of uniform spaces and uniform maps.
In the following we describe the category $\Coarse$ of coarse spaces which is designed to encode the 
large scale structure of metric spaces. Note that not every topological or uniform space comes from a metric space. Similarly, coarse spaces represented by metric spaces only exhaust a small portion of the category of coarse spaces.

  \begin{ddd}A coarse structure 
$\cC$ on a set $X$ is a collection of   subsets of   $X\times X$ whose elements are called coarse entourages.
One requires the following axioms:
\begin{enumerate}
\item The diagonal $\diag(X):=\{(x,x)\mid x\in X\}$  belongs to $\cC$.
\item The set $\cC$ is closed under forming finite unions and taking subsets.
\item The set $\cC$ is closed under the operations flip $$U\mapsto U^{-1}:=\{(y,x)\mid (x,y)\in U\}$$
and composition
$$(U,V)\mapsto U\circ V:=\{(x,y)\in X\times X\mid (\exists z\in X\mid (x,z)\in U\:\&\: (z,y)\in V)\}\ .$$ 
\end{enumerate}
\end{ddd}
A {\em coarse space} is a pair $(X,\cC)$ of a set with a coarse structure $\cC$. If $(X',\cC')$ is a second coarse space and $f:X\to X'$ is a map of underlying sets, then $f$ is called {\em controlled} if  $f\times f$ sends coarse entourages of $X$ to coarse entourages of $X'$.  
We   thus obtain the category  $\Coarse$  of coarse spaces and controlled maps.

\begin{ex}
A distance function $d$ on $X$ gives rise to  the {\em metric coarse structure} $\cC_{d}$ defined as the smallest coarse structure containing the metric entourages \eqref{sbgbsgbvfsfvfvsfdv}.
 
 Assume that $d'$ is a second distance function on $X$. If for every $s$ in $(0,\infty)$ there exists an $r$ in $(0,\infty)$ such that $d(x,y)<s$ implies $d'(x,y)<r$, then we have $\cC_{d}\subseteq \cC_{d'}$ and the identity
 of $X$ is a morphism $(X,\cC_{d})\to (X,\cC_{d'})$ of coarse spaces.  If in the condition above we can interchange the roles of $d$ and $d'$, then $(X,\cC_{d})=(X,\cC_{d'})$. 
\hB \end{ex}

 \begin{ex}
 If $S$ and $T$ are two finite  generating sets of a group $G$, 
 then we have the equality of coarse structures $\cC_{d_{S}}=\cC_{d_{T}}$ on $G$.  Thus by considering the coarse structure on $G$ determined by any choice of such a generating set and studying  geometric properties of $G$ which only depend on the coarse structure we get rid of the dependence on the choice of the generating set.
 \hB\end{ex}
Any subset $Y$ of a coarse space $(X,\cC)$ has an induced coarse structure given by the coarse entourages
of $X$ that are contained in $Y\times Y$.

For any collection of entourages on $X$ we can consider the smallest coarse structure containing the collection.
A coarse structure represented by a metric admits a countable set of generators, namely the family 
 of metric entourages $(U_{n})_{n\in \nat}$. In the case of a path metric space even the one-member family   $(U_{1})$ suffices to generate the coarse structure.
 
 \begin{ex}
 If $((X_{i},\cC_{i}))_{i\in I}$ is a family of coarse spaces, then the {\em free union} $\bigsqcup^{\free}_{i\in I}(X_{i},\cC_{i})$
 is the set $X:=\bigsqcup_{i\in  I} X_{i}$  with the coarse structure generated by the entourages $\bigcup_{i\in I}U_{i}$ for all $(U_{i})_{i\in I}$ in $\prod_{i\in I} \cC_{i}$. If $I$ is infinite, then except for degenerate cases  the coarse structure on the free union is not
 countably generated and therefore does not come from a metric.  \hB\end{ex} 
 
 Note that the category of coarse spaces $\Coarse$ still captures the full information about the underlying sets, i.e. there is a forgetful functor $\Coarse\to \Set$. This functor has left- and right adjoints
 which send a set $S$ to the coarse space $S_{min}$ with coarse structure generated by $\diag(S)$,  or $S_{max}$ with the 
  maximal coarse structure given by the power set of $S\times S$.  Honest large scale geometry starts with the introduction of the notion of a {\em coarse equivalence} which will be explained in \cref{w3rijgowregwregwre9}.

 \section{Bornological coarse spaces}
 
For most applications of large scale geometry the coarse structure needs to be complemented by a notion of 
local finiteness. To this end one introduces the notion of a bornology. \begin{ddd} A bornology  $\cB$ on a set $X$ is a collection of subsets of $X$ called the bounded subsets. One requires the following axioms:
\begin{enumerate}
\item Every finite subset belongs to $\cB$.
\item $\cB$ is closed under forming finite unions and taking subsets.
\end{enumerate} \end{ddd}
A {\em bornological space} is a pair $(X,\cB)$ of a set $X$ with a bornology $\cB$. If $(X',\cB')$ is a second bornological space, then a map $X\to X'$ is called {\em proper} if preimages  of bounded sets   are bounded. It is called {\em bornological} if 
images  of bounded sets   are bounded.

Any set $X$ has a minimal bornology consisting of all finite subsets and a maximal bornology where all subsets are bounded.

A distance function on $X$ determines the {\em metric bornology} $\cB_{d}$  defined as the smallest bornology
containing the metric balls $B(x,r) $ from \eqref{vwewervijovwvwervwerv} for all $x$ in $X$ and $r$ in $(0,\infty)$.
Any subset $Y$ of a bornological space has an induced bornology consisting of the bounded subsets of $X$ which are contained in $Y$.

We consider the bornological space $(X,\cB)$.  \begin{ddd}\label{wethpkjooperthwgrewg} A subset $Y$ of $X$ is called {\em locally finite} if the induced bornology on $Y$ is the minimal one.  \end{ddd}
The collection of locally finite subsets of $X$ again forms a bornology $\cB^{\perp}$.
If $f:(X,\cB)\to (X',\cB')$ is proper, then $f:(X,\cB^{\perp})\to (X',\cB^{',\perp})$ is bornological.

    A  coarse structure $\cC$ and a bornology $\cB$ on the same set $X$ are said to be {\em compatible}  if for every bounded subset $B$ of $X$ and coarse entourage $U$ of $X$ the {\em thickening}  \begin{equation}\label{wetgwetgwrfrewfwrfrfw}U[B]:=\{y\in X\mid (\exists x\in B\mid (y,x)\in U)\}
\end{equation}  is again bounded.  
 
 \begin{ddd}[{\cite[Def. 2.7]{buen}}] A bornological coarse space is a triple $(X,\cC,\cB)$ of a set $X$ with a coarse structure $\cC$ and a  compatible bornology $\cB$.
 A morphism between bornological coarse spaces $f:(X,\cC,\cB)\to (X',\cC',\cB')$ is map of sets
 $f:X\to X'$ which is controlled and proper.
 \end{ddd}
In this way we get the category $\BC$ of bornological coarse spaces.

The metric coarse structure and bornology associated to a distance function $d$ on $X$ are compatible. 
We will denote the associated bornological coarse space by $X_{d}$.

\begin{ex}
For a set $X$ we have the bornological coarse spaces $X_{min,min}$,   $X_{min,max}$ and $X_{max,max}$,
where the first subscript indicates the coarse structure and the second the bornology.
If $X$ is infinite, then the maximal coarse structure and the minimal bornology on $X$ are not compatible. \hB \end{ex} 
\begin{ex} The {\em free union} $\bigsqcup_{i\in I}^{\free}(X_{i},\cC_{i},\cB_{i})$ of a family of bornological coarse spaces $((X_{i},\cC_{i},\cB_{i}))_{i\in I}$ is the free union of the underlying coarse spaces with the bornology generated by $\bigcup_{i\in I}\cB_{i}$. For example, $X_{min,min}$ is the free union of the family of one-point spaces  $(\{x\})_{x\in X}$, while
$X_{min,max}$ is the coproduct. \hB \end{ex}

The category $\BC$ has a {\em symmetric monoidal structure} $\otimes $ defined such that
 $(X,\cC,\cB)\otimes (X',\cC',\cB')$ is given by
 $(X\times X',\cC'',\cB'')$, where
 $\cC''$ is  generated by all entourages $U\times U'$ with $U$ in $\cC$ and $U'$ in $\cC'$, and $\cB''$
is generated by $B\times B'$ with $B$ in $\cB$ and $B'$ in $\cB'$. 
Note the forgetful functor $\BC\to \Coarse$ is symmetric monoidal if we equip the target with the cartesian structure.  On the other hand, $\otimes$ differs from the cartesian structure on $\BC$.

In coarse geometry group actions play an important role. 
If $G$ is a group, then  we can consider the symmetric monoidal category  $\Fun(BG,\BC)$ of  bornological coarse spaces  with a $G$-action by automorphisms. It contains the 
full symmetric monoidal subcategory
 $G\BC$ of {\em $G$-bornological coarse spaces}  $(X,\cC,\cB)$ 
 characterized by the condition that  the subset of $G$-invariant entourages $\cC^{G}$ is cofinal in $\cC$
 with respect to the inclusion relation, i.e., every element of $\cC$ is contained in a $G$-invariant one.
 In this case we say that $\cC$ is a {\em $G$-coarse structure}.

\begin{ex} If $G$ acts isometrically on a metric space $(X,d)$, then the associated  bornological coarse space
 $X_{d}$   belongs to $G\BC$.  \hB
 \end{ex}
 
 \begin{ex} 
If we consider the action of $\Z$ on the metric space $(\R,d_{eu})$ given by $(n,x)\mapsto 2^{n}x$, then  $ \R_{d_{eu}} $ with this action belongs to $\Fun(BG,\BC)$, but not to
 $G\BC$. \hB
 \end{ex}
 
 \begin{ex} 
 The group $G$ itself has a {\em canonical $G$-coarse structure} $\cC_{can}$ generated by
  the family of invariant entourages $(\{(gh,gk)\mid g\in G\})_{(h,k)\in G\times G}$.
 Together with the minimal bornology we get an object $G_{can,min}$ in $G\BC$.
 If $G$ admits a finite   generating set $S$, then $G_{can,min}=G_{d_{S}}$.  \mbox{}\hB
  \end{ex}

 If $(X,\cC)$ is a coarse space, then there is a {\em minimal compatible  bornology}  on $X$.
 The classical literature  only considers bornological coarse spaces whose bornology 
 is the minimal one compatible with the coarse structure. But following examples show that it is useful to decouple the choice of the bornology from the coarse structure.
 
 \begin{ex}
The {\em orbit category} $G\Orb$ of a group $G$ is the category of transitive $G$-sets and equivariant maps. We have a functor
$$i:G\Orb\to G\BC , \quad S\mapsto S_{min,max}\ .$$  
If $S$ is infinite, then the maximal bornology is different from the minimal one compatible with the coarse structure.
This functor is the starting point for the application of coarse geometry to the study of assembly maps \eqref{sfbddopjkpowegfs}. \hB\end{ex} 

\begin{ex}\label{wtrhgwrtgreregwefeer} 
Let $(M,g)$ be a complete Riemannian manifold with an action of a discrete group $G$ by isometries. The Riemannian metric $ g$
allows to measure the length of curves and induces a Riemannian distance function $d_{ g}$, and one usually considers the $G$-bornological coarse space $M_{d_{g}}=(M,\cC_{d_{g}},\cB_{d_{g}})$. But  in   index theory of Dirac operators it will be useful to work with a larger bornology $\cB_{s_{g}}$.
Let $s_{g}:M\to \R$ be the scalar curvature function. Then the bornology $\cB_{s_{g}}$ consists of the subsets $B$ of $M$ with the property  that there exists a compact subset $K$ of $M$ such that $\inf_{B\setminus K} s_{g}>0$. One checks that $\cB_{s_{g}}$ is compatible with
the metric coarse structure and we get a bornological coarse space $M_{d_{g},s_{g}}:=(M,\cC_{d_{g}},\cB_{s_{g}})$. The identity of $M$ is a morphism $
M_{d_{g},s_{g}}\to M_{d_{g}}$ in $G\BC$. 
We will explain in  \cref{qrfgqoirejfqofeweqfd} that this example allows to capture the fact that
the coarse index of an invariant $Spin$ Dirac operator is supported away from the set where the scalar curvature is positive. \hB
\end{ex}

\section{Coarse equivalence}\label{w3rijgowregwregwre9}

Two morphisms $f_{0},f_{1}:X\to Y$  between $G$-bornological coarse spaces are said to be {\em close} to each other if $(f_{0}\times f_{1})(\diag(X))$ is a coarse entourage of $Y$. Equivalently we can require that the map
$\{0,1\}_{max,max}\otimes X\to Y$ given by $(i,x)\mapsto f_{i}(x)$ is a morphism in  $G\BC$.
Closeness  is an equivalence relation on morphisms   which is compatible with composition. 

\begin{ddd}
A map $f:X\to Y$ in $G\BC$ is a coarse equivalence if it is invertible up to closeness.
\end{ddd}
In detail this means that there exists a map $h:Y\to X$ such that $h\circ f$ is close to $\id_{X}$ and $f\circ h$ is close to $\id_{Y}$.

A subset $Y$ of a bornological coarse space $X$ is called {\em dense} if there exists a coarse entourage $U$ of $X$ such that $U[Y]=X$, see \eqref{wetgwetgwrfrewfwrfrfw}. In this case the inclusion $Y\to X$ is a coarse equivalence, where we equip $Y$ with the induced bornological coarse structures.
 In the equivariant case this is not always true.  \begin{ex} Consider the $C_{2}$-bornological coarse space $\R_{d_{eu}}$  such that the non-trivial element of $C_{2}$  acts by multiplication by $-1$. Then
 the inclusion $\R_{d_{eu}}\setminus\{0\}\to \R_{d_{eu}}$ is dense and a coarse equivalence of underlying bornological
 coarse spaces, but not a coarse equivalence in $C_{2}\BC$ since there exists no equivariant map of sets 
 $ \R\to \R\setminus\{0\}$ at all. 
 In view of this example one sometimes considers the notion of a {\em weak coarse equivalence} (see e.g. \cite[Def. 2.18]{coarsek}) between $G$-bornological coarse spaces which is map which becomes a coarse equivalence after forgetting the $G$-action.
  \hB\end{ex}

 \begin{ex} \label{3jqrrgwregegw}Let $(M,g)$ be a connected compact Riemannian manifold with fundamental group $G$. Its universal covering $(\tilde M,\tilde g)$ is a complete Riemannian manifold with an isometric $G$-action.  It has an induced $G$-invariant distance function $d_{\tilde g}$. If $\tilde m_{0}$ is any point in $\tilde M$, then the map
 $G\to \tilde M$, $g\mapsto gm_{0}$, induces a coarse equivalence $G_{can,min}\to \tilde M_{d_{\tilde g}}$ in $G\BC$.  \hB  \end{ex}

As a general rule, all concepts of coarse geometry should be invariant under coarse equivalences.
In the following we provide two examples, the set of coarse components and the Higson corona.
Further examples of coarsely invariant  concepts  are  the notions of {\em bounded geometry} \cite{hr}, {\em property $A$} \cite{Yu2000},
  {\em asymptotic dimension} \cite{gromov}, or {\em finite decomposition complexity} \cite{Guentner:2010aa}.

  \begin{ex}
In the following we  describe  the   functor $$\pi_{0}^{\coarse}:G\BC\to G\Set$$
which associates to every $G$-bornological coarse space its $G$-set of coarse components. 
Let $X$ be a  $G$-bornological coarse space. 
\begin{ddd} The $G$-set of  coarse components $\pi_{0}^{\coarse}(X)$   is the set of equivalence classes on $X$ with respect to the 
equivalence relation $R_{\cC}:=\bigcup_{U\in \cC}U$ and induced $G$-action, where $\cC$ denotes the coarse structure of $X$. \end{ddd}
Note that $\pi_{0}^{\coarse}(X)$
 is independent of the bornology.
 A map of $G$-bornological coarse spaces $f:X\to Y$ functorially  induces a map of $G$-sets  $$\pi_{0}^{\coarse}(f):\pi_{0}^{\coarse}(X)\to \pi_{0}^{\coarse}(Y) \ .$$
If  $f_{0}$ is close to $f_{1}$, then we have $\pi_{0}^{\coarse}(f_{0})=
\pi_{0}^{\coarse}(f_{1})$. So the functor $\pi_{0}^{\coarse}$  is coarsely invariant. It even sends weak coarse equivalences to isomorphisms.

For a $G$-set $X$ we have a canonical isomorphism  $X \cong \pi_{0}^{\coarse}(X_{min,min})$. The group $G$ with its canonical coarse structure is coarsely connected, i.e., we have  $\pi_{0}^{\coarse}(G_{can,min})\cong *$.
If the $G$-coarse structure comes from an invariant  distance function $d$ on $X$, then $x$ and $y$ belong to the same coarse component of $X$ if and only if $d(x,y)<\infty$.  
 \hB
\end{ex}
 
 \begin{ex}
 In this example we discuss the  {\em Higson corona functor}
 \begin{equation}\label{rgwregrfrefrwf}\partial: G\BC\to G\Top_{\mathrm{cp}}
\end{equation} 
 which associates to a $G$-bornological coarse space $X$ a compact Hausdorff space $\partial X$ with $G$-action.
Let  $\phi:X\to \C$ be a bounded function. For any entourage $U$  and   subset $Y$ of $X$  we define the $U$-variation of $\phi$ on $Y$ by
$$\Var_{U}(\phi)(Y):=\sup_{(y,y')\in U\cap (Y\times Y)} |\phi(y)-\phi(y')| \ .$$
We then consider the $G$-$C^{*}$-subalgebra $$C_{h}(X):=\{\phi\in C_{b}(X)\mid (\forall U\in \cC\mid \lim_{B\in \cB} \Var_{U}(\phi)(X\setminus B)=0)\}$$ of bounded functions with vanishing variation at infinity. 
We furthermore consider the $G$-$C^{*}$-subalgebra  $C_{0}(X)$ generated by  the functions with  bounded support.
 We finally define the unital $G$-$C^{*}$-algebra \begin{equation}\label{qwefeqeqwfqwdada}C(\partial X):=C_{h}(X)/C_{0}(X)
\end{equation}
 \begin{ddd}[{\cite{hr}}]\label{fwerfwrfrefw5w}
 The Higson corona $\partial X$ of $X$ is the compact Hausdorff space corresponding to $C(\partial X)$ under Gelfand duality with the induced $G$-action.
 \end{ddd}

  If $f:X\to Y$ is a morphism of $G$-bornological coarse spaces, then the canonical homomorphism $f^{*}:C_{b}(Y)\to C_{b}(X)$ sends the subalgebras $C_{h}(Y)$ and $C_{0}(Y)$ to $C_{h}(X)$ and $C_{0}(X)$, respectively and therefore
  induces an equivariant homomorphism of quotients $\bar f^{*}:C(\partial Y)\to C(\partial X)$ and hence a continuous map of coronas $\partial f:\partial X\to \partial Y$. We thus get the functor \eqref{rgwregrfrefrwf}.
  
  If $f_{0},f_{1}:X\to Y$ are close to each other and $\phi$ is in $C_{h}(Y)$, then $f_{0}^{*}\phi-f_{1}^{*}\phi\in C_{0}(X)$. This implies that $\partial f_{0}=\partial f_{1}$. Therefore the corona functor  $\partial$ is a coarsely invariant.  
 
 Applying properties of compact Hausdorff spaces  to the Higson corona we get a variety of coarsely invariant concepts.
 
In \cref{qrfgqoirejfqofeweqfd} we will explain that the  topological $K$-theory $K(\partial X)$ of the Higson corona 
  pairs interestingly  with the  coarse $K$-homology
 $K\cX(X)$ of $X$.
  \hB
 \end{ex}

%The corona depends on the bornology. We consider the metric space $\R_{d}$.  

%For example $\partial X_{min,max}$ is the Stone-\v{C}ech compactification of the discrete space $X$, while $\partial X_{min,min} $

\section{Coarse homology theories}

In the following we describe an axiomatization of   equivariant coarse homology theories.
 We consider a functor $E:G\BC\to \bM$ whose target $\bM$ is a cocomplete stable $\infty$-category, e.g. the category of spectra $\Sp$ or the derived category of abelian groups $D(\Z)$.
\begin{rem} Following \cite[Sec. 1]{HA}
an $\infty$-category $\bM$ is called stable, if it is pointed and admits finite limits and colimits, and if push-out squares in $\bM$ are the same as pull-back squares.   %Define the $\infty$-category $\Spc$ of spaces as the Dwyer-Kan localization of the category of topological spaces at the weak homotopy equivalences. 
 The $\infty$-category of spectra $\Sp$ is the universal presentable stable $\infty$-category generated by an object $S$ called the sphere spectrum. If $E$ is a spectrum, then we  define its  $\Z$-graded homotopy groups by
 $$\pi_{n}E:=[\Sigma^{n}S,E]\ ,$$ where $[-,-]$ denote the  group of maps in the homotopy
 category of $\Sp$. In particular, a morphism between spectra is an equivalence if and only if it induces an isomorphism between the homotopy groups.
 
   For any stable $\infty$-category $\bM$  and objects $A,B$  we have  a mapping spectrum $\map_{\bM}(A,B)$ in $\Sp$.
 %Similarly, the $\infty$-category $D(\Z)$ is the universal presentable stable $\Z$-linear $\infty$-category generated by an object, also denoted by  $\Z$.

 The stable $\infty$-category $D(\Z)$   is  the Dwyer-Kan localization of the category of chain complexes of abelian groups at the quasi-isomorphisms. It is the universal  presentable stable $\Z$-linear $\infty$-category generated by $\Z$ considered as a chain complex in the natural way. If $A$ is any chain complex, then  we have  an isomorphism
 $$\pi_{n}\map_{D(\Z)}(\Z,A)\cong H_{n}(A)\ .$$
 Again, a morphism in $D(\Z)$ is an equivalence if and only if it induces an isomorphism of homology groups.
 %In the case of spectra and for $D=S$ the long exact sequence from above is the usual long exact sequence of homotopy groups associated to a fibre sequence of spectra. In $D(\Z)$ one can represent
%the fibre sequence by a short exact sequence of chain complexes. If we take $D=\Z$, then the long exact sequence is the usual long exact sequence of homology groups of chain complexes known
%from homological algebra.
\hB 
\end{rem}
 
 %For any  fibre sequence $A\to B\to C$      and object $D$ in $\bM$ we get a long exact sequence of abelian groups
%$$\dots \to \pi_{*}\map_{\bM}(D,A)\to \pi_{*}\map_{\bM}(D,B)\to \pi_{*}\map_{\bM}(D,C)\to \pi_{*-1}\map_{\bM}(D,A)\to \dots\ . $$

 %Using the constructions from above one can derive the classical group valued functors and the corresponding long exact sequences in a canonical way.
 Classically, coarse homology theories are defined as $\Z$-graded abelian group valued functors. In this case the
boundary operators for Mayer-Vietoris sequences have to be constructed as an {\em additional datum}.  
Working with functors with values in a stable $\infty$-category turns Mayer-Vietoris into a {\em property}    of the functor $E:G\BC\to \bM$, see  \cref{weijotgtggd9}.\ref{gjieogerferwf} below for a precise formulation.

\begin{ddd}[{\cite[Def. 3.10]{equicoarse}}]\label{weijotgtggd9}
$E$ is an equivariant coarse homology theory if
\begin{enumerate}
\item $E$ is coarsely invariant.
\item\label{gjieogerferwf} $E$ is excisive.
\item $E$ annihilates flasques.
\item $E$ is $u$-continuous.
\end{enumerate}
\end{ddd}
In the following we explain  the meaning of these conditions.

{\em Coarse invariance} means that $E$ sends coarse equivalences to equivalences.
This could equivalently be phrased as the condition that the projection onto $X$ induces an equivalence 
$$\pr:E(\{0,1\}_{min,max}\otimes X)\stackrel{\simeq}{\to} E(X)$$
for every $G$-bornological coarse space $X$. 
Coarse invariance is the main condition which ensures that associating $E(X)$ to $X$ is a coarsely invariant concept. There are interesting examples of equivariant coarse homology theories  which even send weak coarse equivalences to equivalences, e.g. \cite[Thm. 8.7]{coarsek} and the examples from \cite[Cor. 5.3.13]{unik}.

{\em Excisiveness}  is the condition which makes  $E$ homological. It means that
 $E(X)$ can be determined by glueing the values of $E$ on   pieces of $X$. For a precise formulation
we introduce the notion of a {\em complementary pair} $(Z,\cY)$ in $X$. Here $Z$ is an $G$-invariant subset of $X$ and $\cY$
is a family $(Y_{i})_{i\in I}$ of $G$-invariant subsets indexed by a filtered poset $I$ which is big in the sense that
for any $i$ in $I$ and invariant coarse entourage $U$ of $X$ there exists $i'$ in $I$ such that $U[Y_{i}]\subseteq Y_{i'}$, see \eqref{wetgwetgwrfrewfwrfrfw}.
In addition we require that there exists $i$ in $I$ with $Y_{i}\cup Z=X$.
We define $E(\cY):=\colim_{i\in I} E(Y_{i})$ and consider the big family $Z\cap \cY:=(Z\cap Y_{i})_{i\in I}$ on $Z$. Then $E$ is called excisive if
$$ \xymatrix{E(Z\cap \cY)\ar[r]\ar[d] &E(Z) \ar[d] \\E(\cY) \ar[r] &E(X) }$$ 
is a push-out square in $\bM$ for every $G$-bornological coarse space $X$ and complementary pair $(Z,\cY)$.
If $(Z,Y)$ is an invariant covering of $X$, then one considers the complementary pair $(Z,\{Y\})$ with the big family
$\{Y\}:=(U[Y])_{U\in \cC^{G}}$ of all invariant thickenings of $Y$. We   work with the big family $\{Y\}$ instead of $Y$  since
the intersection $Z\cap \{Y\}$ is a coarsely invariant concept in contrast to $Z\cap Y$.
But in many examples the canonical map $E(Y)\to E(\{Y\})$ is an equivalence since the inclusions $Y\to U[Y]$ are coarse equivalences.
 
A $G$-bornological coarse space $X$ is called {\em flasque} if it admits a selfmap $f:X\to X$ such that
$f$ is close to $\id_{X}$, $f$ is non-expanding in the sense that for any coarse entourage $U$ of $X$ the union $\bigcup_{n\in \nat} (f^{n}\times f^{n})(U)$
is again a coarse entourage of $X$, and $f$ shifts $X$ to infinity in the sense that for every bounded subset $B$ of $X$ there exists an $n$ in $\nat$ such that $B\cap f^{n}(X)=\emptyset$.
A typical flasque space is $[0,\infty)  \otimes X$ with $f$ defined by $f(n,x):=(n+1,x)$, where $[0,\infty) $ has the bornological coarse structure induced from the inclusion into $\R_{d_{eu}}$.
We say that $E$ {\em annihilates flasques} if $E(X)\simeq 0$ for every flasque bornological coarse space.
 Annihilation of flasques reflects a version of local finiteness of $E$.
 
 If $X$ is a $G$-bornological space with coarse structure $\cC$, then for any $U$ in $\cC^{G}$ we can consider
 the $G$-bornological coarse space $X_{U}$ with the smaller $G$-coarse structure generated by $U$. The identity of $X$ induces maps $X_{U}\to X$ and $X_{U}\to X_{U'}$ for all $U'$ in $\cC^{G}$ with  $U\subseteq U'$. We say that $E$ is {\em $u$-continuous} 
 if the canonical map is an equivalence  $$\colim_{U\in \cC^{G}} E(X_{U})\stackrel{\simeq}{\to} E(X)$$
 for every $G$-bornological coarse space $X$.
 
 If we fix any compact object $M$ of $\bM$, then we can form a coarsely invariant $\Z$-graded group-valued functor  $X\mapsto \pi_{*} \map_{\bM}(M,E(X))$ which is $u$-continuous and  vanishes on flasques. The excisiveness of $E$ gives rise to long exact Mayer-Vietoris sequences for complentary pairs. 
 
 The classical examples of coarse homology theories are usually defined as   group valued functors
 on certain subcategories of $G\BC$. But most of them can be extended to all of $G\BC$ and admit a
model in the sense defined above.  We refer to \cite{equicoarse} for {\em coarse equivariant ordinary homology} (see also \cref{wioegogfre9}), to  \cite{equicoarse}, \cite{buci} for  {\em coarse algebraic $K$ with coefficients in an additive category}, to \cite{Bunke:aa} and \cite{unik}
for {\em coarse algebraic $K$-theory of spaces} and {\em coarse algebraic $K$ with coefficients in a left-exact $\infty$-category}, to \cite{coarsek} for  {\em equivariant coarse topological $K$-theory} (see also \cref{qrfgqoirejfqofeweqfd}), and to \cite{Caputi_2020} for {\em equivariant coarse cyclic and Hochschild homology}. In the non-equivariant case every locally finite homology theory
admits a {\em coarsification} \cite[Sec. 5.5]{roe_lectures_coarse_geometry}, \cite[Sec. 7]{buen}, \cite{ass}.

One can derive the  notion of an {\em equivariant coarse cohomology} by dualizing the axioms in \cref{weijotgtggd9}, see \cite{Bunke:2017aa}. For alternative sets of axioms, mainly for group-valued functors, see e.g.
 \cite{MR1834777}, \cite{wulff_axioms}.

\begin{ex}
In order to demonstrate the usage of the axioms we do some calculations. Let $E$ be an equivariant  coarse homology theory.
We first show that  for any $X$ in $G\BC$ we have an equivalence
 \begin{equation}\label{sdfbposkfodpbsfdvsvsfdvsfdv}E(\R_{d_{eu}}^{n}\otimes X)\simeq \Sigma^{n} E(X)\ .
\end{equation} 
 This can be shown by induction on $n$. The basic step uses the complementary pair $((-\infty,0]\times X,([-n,\infty)\times X)_{n\in \nat})$ on $\R_{d_{eu}}\otimes X$. By excision we have  a push-out square  $$
  \xymatrix{\colim_{n\in \nat}E([-n,0]\otimes X)\ar[r]\ar[d] &E((-\infty,0]\times X) \ar[d] \\ \colim_{n\in \nat } E([-n,\infty)\times X)\ar[r] & E(\R_{d_{eu}}\otimes X)}  \ .$$
 By coarse invariance the upper left corner is equivalent to $E(X)$, and the lower left and upper right corners are zero since $E$ vanishes on flasques. Hence this square is equivalent to $$  \xymatrix{ E( X)\ar[r]\ar[d] &0 \ar[d] \\ 0\ar[r] &  E(\R_{d_{eu}}\otimes X)} \ .$$ \hB
 \end{ex}
  
  \begin{ex} Let $G$ be trivial.
We consider the subspace $Sq:=\{n^{2}\mid n\in \nat\}$ of $\R_{d_{eu}}$ of square integers in $\R$.
For every entourage $U$ we have a decomposition of  $Sq_{U}$ into a bounded connected subset and
an infinite discrete subset. Let us assume that $E$ is additive in the sense that
$E(X_{min,min})\cong \prod_{X}E(*)$  (see \cite[Def. 6.4]{buen}) for any set $X$. All examples mentioned above have this additional property. 
We then have 
$$E(Sq_{U})\simeq E(*)\oplus\prod_{i=n}^{\infty}E(*)$$ for a suitable $n$. For $n$ in $\nat$ we consider the map
$$ E(*)\oplus\prod_{i=n}^{\infty}E(*)\to E(*)\oplus \prod_{i=n+1}^{\infty} E(*)\ ,  \quad (a,(x_{i})_{i\in \{n,\dots,\infty\}})\mapsto  (a+x_n,(x_{i})_{i\in \{n+1,\dots,\infty\}})\ .$$
Then using $u$-continuity of $E$ we get
$$E(Sq)\simeq \colim_{n\in \nat } \left( E(*)\oplus\prod_{i=n}^{\infty}E(*)\right)\simeq   E(*)\oplus \frac{\prod_{\nat}E(*)}{\bigoplus_{\nat}E(*)}\ .$$ \hB
\end{ex}

 \section{Equivariant coarse ordinary homology}\label{wioegogfre9}
 
 A version of ordinary coarse (co)homology was first defined by Roe \cite{MR1147350}. For further constructions see \cite{block_weinberger_large_scale}, \cite{Schmidt_1999}, \cite{Hartmann_2020}.
 Following \cite[Sec. 7]{equicoarse}, in this section we sketch the construction of the equivariant coarse homology functor
$$H\cX^{G}:G\BC\to D(\Z) $$  in the sense of \cref{weijotgtggd9}.
%where $D(\Z)$  is the derived  (stable $\infty$-)category of $\Z$-modules.
Let $X$ be a $G$-bornological coarse space. Let $U$ be a coarse entourage of $X$ and $B$ be a bounded subset. A point $(x_{0},\dots,x_{i})$ in $X^{n+1}$ is called $U$-controlled if $(x_{i},x_{i'})\in U$ for all   $i,i'$ in $\{0,\dots,n\}$. The point meets $B$ if $x_{i}\in B$ for some $i$ in $ \{0,\dots,n\}$.

We consider the group $C\cX_{n,U}^{G}(X)$ of all $G$-invariant $U$-controlled and locally finite  functions
 $\phi:X^{n+1}\to \Z$. The latter two properties require that the support of $\phi$ consists of $U$-controlled points and that every bounded subset $B$ meets the support of $\phi$ in a finite subset. We define the differential
 $\partial:C\cX_{n+1,U}^{G}(X)\to C\cX_{n,U}^{G}(X)$ by the usual formula
 $$\partial \phi(x_{0},\dots,x_{n}):=\sum_{i=0}^{n+1} (-1)^{i} \sum_{x\in X} \phi(x_{0},\dots,x_{i-1},x,x_{i},\dots,x_{n})\ .$$
 The conditions on $\phi$ ensure  that the sum has only finitely many non-zero terms.
 For every invariant entourage $U$ we get a chain complex $C\cX^{G}_{U}(X)$ which we consider as an object of $D(\Z)$. If $X\to X'$ is a proper map such that $(f\times f)(U)\subseteq U'$, then we get a map of chain complexes 
 $$f_{*}:C\cX^{G}_{U}(X) \to C\cX^{G}_{U'}(X')\ , \quad (f_{*}\phi)(x_{0}',\dots,x'_{n}):=\sum \phi(x_{0},\dots,x_{n})$$ 
 where the sum is taken over the fibre of the map $f^{n+1}:X^{n+1}\to X^{\prime,n+1}$. It again has finitely many non-zero terms.
We then define $$H\cX^{G}(X):=\colim_{U\in \cC^{G}} C\cX^{G}_{U}(X)$$ in $D(\Z)$.
The construction is functorial in $X$.

\begin{ddd}
The functor $H\cX^{G}:G\BC\to D(\Z)$ is called the equivariant coarse ordinary   homology theory.
\end{ddd}

It satisfies the axioms from \cref{weijotgtggd9}. %By taking homology we  get  classical ordinary coarse homology theory functor  $H\cX^{G}_{*}(-)$ together with its natural Mayer-Vietoris sequences

 \begin{ex}
 In the case of the trivial group we omit $G$ from the notation. 
 By an explicit calculation we have $$H\cX_{k}(*)= \left\{\begin{array}{cc} \Z&k=0\\0 &\mbox{else}  \end{array} \right.\ .$$
By specializing  \eqref{sdfbposkfodpbsfdvsvsfdvsfdv} we  get  $$H\cX_{k}(\R^{n}_{d_{eu}})\cong \left\{\begin{array}{cc} \Z& k=n\\0& \mbox{else}  \end{array} \right.$$ \hB
\end{ex}

\begin{ex} For non-trivial groups equivariant coarse homology is related with group homology via
$$H\cX^{G}(G_{can,min}\otimes S_{min,max})\simeq H(G,\Z[S])\ ,$$
where $S$ is a $G$-set, $\Z[S]$ in $G\Mod(\Z)$ is the associated $G$-module, and $H(G,-)$ is the group homology functor
$G\Mod(\Z)\to D(\Z)$  \cite[Prop. 3.8]{engel_loc}, \cite[Prop. 7.5]{equicoarse}. \hB
\end{ex} 
\section{Equivariant coarse topological $K$-homology}\label{qrfgqoirejfqofeweqfd}
Group-valued coarse topological $K$-homology for proper metric spaces  is usually defined in terms of the $K$-theory of Roe algebras \cite{roe_lectures_coarse_geometry},  \cite{hr}, \cite{willett_yu_book}.
In this section, following \cite{coarsek} we sketch the construction of a spectrum-valued  equivariant coarse $K$-homology theory in the sense of  \cref{weijotgtggd9} using Roe $C^{*}$-categories.

Let $X$ be a $G$-bornological coarse space. An {\em $X$-controlled Hilbert space} is a triple $(H,\rho,p )$ of a Hilbert space $H$, a unitary representation $\rho=(\rho_{g})_{g\in G}:G\to U(H)$, and a mutually orthogonal family of  finite dimensional  orthogonal projections $p :=(p_{x})_{x\in X}$ on $H$ such that $\sum_{x\in X} p_{x}=\id_{H}$ strongly,
   the set $\{x\in X\mid p_{x}\not=0\}$ is locally finite (see \cref{wethpkjooperthwgrewg}), and
such that $\rho_{g}p_{x} \rho_{g}^{-1}=p_{gx}$ for all $g$ in $G$ and $x$ in $X$.
 
A controlled morphism between  $X$-controlled Hilbert spaces $A:(H,\rho,p )\to (H',\rho',p')$
is a bounded operator $A:H\to H'$ such that $\rho'_{g}A=A\rho_{g}$ for all $g$ in $G$, and such that the set
$\{(x,y)\in X\times X\mid p'_{x}Ap_{y}\not=0\}$ is a coarse entourage of $X$.
We let $C^{*}((H,\rho,p),(H',\rho',p'))$ denote the closure in $B(H,H')$ of the set of  controlled morphisms from $(H,\rho,p)$ to $(H',\rho',p')$. 

We obtain a $C^{*}$-category $\bV^{G}(X)$ \cite{ghr}, \cite{cank} whose objects are the   $X$-controlled Hilbert spaces, and whose morphism spaces are the spaces $ C^{*}((H,\rho,p),(H',\rho',p'))$ with the obvious involution given by taking the adjoint operator.  

If $X$ itself is locally finite and $\dim(p_{x})=1$ for all $x$ in $X$, then 
the endomorphism $C^{*}$-algebra
$$C^{*}(H,\rho,p):=\End_{\bV^{G}(X)} ((H,\rho,p))$$
is called  the {\em equivariant uniform Roe algebra}  associated to $X$.

 If $f:X\to X'$ is a morphism in $G\BC$, then we get a functor
$$\bV^{G}(f): \bV^{G}(X)\to \bV^{G}(X')$$ which sends $(H,\rho,p)$ to $(H,\rho,f_{*}p)$ with
$(f_{*}p)_{x'}:=\sum_{x\in f^{-1}(x')} p_{x}$ and the morphism $A$ in $ C^{*}((H,\rho,p),(H',\rho',p'))$ to the same operator
$A$, now  considered as an element in $ C^{*}((H,\rho,f_{*}p),(H',\rho',f_{*}p'))$.

We get a functor
$$\bV^{G}:G\BC\to \Ccat\ .$$
 We now use the   spectrum-valued $K$-theory functor for $C^{*}$-categories \begin{equation}\label{fvwovjeorivjevvsdvsdfvsfv}K:\Ccat\to \Sp\ ,
\end{equation} see  \cite{joachimcat}, \cite{cank} for details. \begin{ddd} We define 
  the equivariant   coarse topological $K$-homology functor by
 \begin{equation}\label{sfdbsdvsfdvfdvsdrrf}K\cX^{G}:=K\circ \bV^{G}:G\BC\to \Sp\ .
\end{equation} \end{ddd}
 By \cite[Thm. 7.3]{coarsek} this functor is an equivariant coarse homology theory  in the sense of \cref{weijotgtggd9}. On  bornological coarse spaces associated to proper metric spaces with isometric group actions  the corresponding group valued functor coincides with the classical
 construction of equivariant coarse $K$-homology \cite{roe_lectures_coarse_geometry}, see \cite[Th. 6.1]{indexclass}.

 \begin{ex}
By an explicit calculation we get 
 $$K\cX(*)\simeq KU\ ,$$
 where $KU$ denotes the complex $K$-theory spectrum. By  specializing  \eqref{sdfbposkfodpbsfdvsvsfdvsfdv} 
 we get \begin{equation}\label{qregjqnregioergergrwe}K\cX(\R^{n}_{d_{eu}})\simeq \Sigma^{n}KU\ .
\end{equation}
 \hB\end{ex}
 
 \begin{ex}\label{werguhw9egweferwfw}
 The functor \begin{equation}\label{rwerffwfqfqewfwf}K^{G}:G\Orb\to \Sp\ , \quad  S\mapsto K\cX^{G}(S_{can,min}\otimes G_{can,min})
\end{equation} 
 is the functor constructed in \cite{davis_lueck}.
 Its values are given by  $K^{G}(G/H)\simeq K(C^{*}_{r}(H))$  \cite[Prop. 9.2.3]{coarsek}, where $C_{r}^{*}(H)$ is the reduced group $C^{*}$-algebra of the subgroup $H$ of $G$.
 The functor $K^{G}$
  features the left-hand side of the Baum-Connes/Davis-Lück assembly map \begin{equation}\label{sfbddopjkpowegfs}\colim_{G_{\Fin}\Orb} K^{G}\to K^{G}(*)\ ,
\end{equation}  (see \cite{kranz}, \cite{bel-paschke} for the comparison of the two versions of assembly maps),
where
$G_{\Fin}\Orb$ is the full subcategory of $G\Orb$ of transitive $G$-sets with finite stabilizers.
 Coarse geometry can be used to show that this assembly map is split injective. For example it is known by \cite{skandalis_tu_yu} that   the assembly map \eqref{sfbddopjkpowegfs} is split injective if $G_{can,min}$ admits a coarse embedding into a  Hilbert space. In \cite{desc} a completely different argument  (axiomatizing \cite{Guentner:2010aa}) applies coarse geometry to show split-injectivity of \eqref{sfbddopjkpowegfs} and its versions for   functors on the orbit category derived from other coarse homology theories as in \eqref{rwerffwfqfqewfwf}  under finite decomposition complexity assumptions on $G$. \hB
 \end{ex}

 \begin{ex}Combining \cref{3jqrrgwregegw} and \cref{werguhw9egweferwfw} for $S=*$ we get an equivalence
 $$K\cX^{G}(\tilde M_{d_{\tilde g}})\simeq K(C^{*}_{r}(G))\ .$$ \hB
   \end{ex}

 \begin{ex}
 The equivariant coarse $K$-homology naturally captures the index of equivariant Dirac type operators on complete Riemannian manifolds $(M,g)$ with a proper action of $G$ by isometries. Assume that $M$ admits an
 equivariant spin structure and let $\Dirac$  denote the associated Dirac operator.  As observed in \cite{Roe:2012fk}, \cite{indexclass} it then has a well-defined index class
 $$\ind\cX(\Dirac)\:\:\mbox{in}\:\: K\cX^{G}_{-\dim(M)}(M_{d_{g},s_{g}})\ .$$
 The appearance of the bornology $\cB_{s_{g}}$ (see \cref{wtrhgwrtgreregwefeer}) reflects the fact that the index class is essentially supported away from the subset of $M$ where the scalar curvature is positive.
 
 If $G$ is trivial and  $s_{g}$ admits a uniform positive  lower bound outside  of a compact subset, then $\Dirac$  is Fredholm, and its Fredholm index $\ind(\Dirac)$ in $K_{-\dim(M)}(\C)$ can be expressed using the functoriality of the coarse $K$-homology
 as $\ind(\Dirac)= p_{*}\ind\cX(\Dirac)$ in $K\cX_{-\dim(M)}(*)\cong \pi_{-\dim(M)}KU $, where $p:M_{d_{g},s_{g}}\to *$ is the projection. \hB \end{ex}
  \begin{ex}
  Recall the Higson corona  defined in  \cref{fwerfwrfrefw5w}.  For simplicity we consider the case of the trivial group. 
  \begin{prop} There exists a binatural pairing
$$K(\partial X)\otimes K\cX(X)\to \Sigma KU\ .$$\end{prop}
In order to construct this pairing 
we let $ \Hilb(\C)/\Hilb_{c}(\C)$ be the Calkin $C^{*}$-category whose objects are Hilbert spaces, and whose morphisms are the quotients of bounded by compact operators. Following \cite{MR2661442}
there is a functor
\begin{equation}\label{dfsvsjevdfssvfdv}C(\partial X)\otimes \bV(X)\to  \Hilb(\C)/\Hilb_{c}(\C)
\end{equation} (the domain is a tensor product of $C^{*}$-categories \cite[Sec. 7]{KKG}) which sends
the object $(H  ,p)$ to $H$ and the morphism $[\phi]\otimes A: (H , p)\to (H'  ,p')$ to $[ Am(\phi) ]:H\to H'$ in 
$ \Hilb(\C)/\Hilb_{c}(\C)$. Here $\phi$ in $C_{h}(X)$ denotes a representative of the class $[\phi]$ in the quotient \eqref{qwefeqeqwfqwdada},  $m(\phi)=\sum_{x\in X} \phi(x)p_{x}$ is strongly convergent in $B(H)$,
and $[ Am(\phi) ]$ denotes the class of  $Am(\phi):H\to H'$ in the Calkin category.
  Using $K(\partial X):=K(C(\partial X))$, the symmetric monoidal structure of the
$K$-theory functor \eqref{fvwovjeorivjevvsdvsdfvsfv} and \eqref{sfdbsdvsfdvfdvsdrrf} for the first map, and \eqref{dfsvsjevdfssvfdv} for the second, the pairing is given by the composition
$$K(\partial X)\otimes K\cX(X)\to K(C(\partial X)\otimes \bV(X))\to  K( \Hilb(\C)/\Hilb_{c}(\C))\simeq \Sigma KU\ .$$
The pairing of  coarse index  classes  of Dirac type operators on complete Riemannian manifolds
with $K$-theory classes  from the Higson corona can be expressed in terms of Callias-type operators 
 \cite{Bunke_1995}.
 In order to exhibit a non-trivial example
 we consider the natural map $p:\partial \R^{n}_{d_{eu}}\to S^{n-1}$. If we choose a generator $o_{S^{n-1}}$ of the reduced $K$-group   $\tilde K^{n-1}(S^{n-1})\cong \Z$, then   the pairing
 $$p^{*}o_{S^{n-1}}\otimes-:K\cX_{n}(\R^{n}_{d_{eu}})\to  \pi_{1}\Sigma KU\cong \Z$$ is an isomorphism reflecting the equivalence \eqref{qregjqnregioergergrwe}. 
  \hB
\end{ex}

\section{Summary}

In this survey we introduced the category of $G$-bornological coarse spaces as a basic framework for coarse geometry.  We explained the concepts of   a coarse equivalence and of  an equivariant coarse homology theory. We gave 
 complete constructions of equivariant coarse ordinary homology and equivariant coarse $K$-homology and indicated some basic calculations.

\paragraph*{Keywords:}  coarse space, bornological coarse space, coarse equivalence, coarse homology theory, Higson corona, coarse ordinary homology, controlled Hilbert space, Roe algebra,  Roe category, coarse $K$-homology, canonical coarse structure on a group, coarse index class, assembly map

%\bibliographystyle{plain}
%\bibliography{ency}

\end{document}